\documentclass[a4paper]{jpconf}
\usepackage[square,sort,compress,numbers]{natbib} % necessary for this bibliographystyle, ensures proper sorting and ``References'' headline
\usepackage{amssymb,amsmath, amsfonts, amsthm, iopams}
\usepackage[active]{srcltx} % For inverse serarch

\newtheorem{definition}{Definition}

\newtheorem{theorem}[definition]{Theorem}
\newtheorem{lemma}[definition]{Lemma}

\newtheorem{proposition}[definition]{Proposition}

\newcommand{\G}{\varGamma}

\renewcommand{\L}{\varLambda}

\newcommand{\ti}[1]{\widetilde{#1}}
\newcommand{\set}[1]{\left\lbrace #1\right\rbrace}

\newcommand{\ts}{\hspace{0.5pt}}

\newcommand{\RR}{\mathbb{R}}
\newcommand{\QQ}{\mathbb{Q}}
\newcommand{\ZZ}{\mathbb{Z}}
\newcommand{\NN}{\mathbb{N}}
\newcommand{\HH}{\mathbb{H}}
\newcommand{\II}{\mathbb{I}}
\newcommand{\oo}{{\scriptstyle \mathcal{O}}}

\renewcommand{\and}{\quad \text{and} \quad}

\newcommand{\ii}{\mathrm{i}}
\newcommand{\jj}{\mathrm{j}}
\newcommand{\kk}{\mathrm{k}}

\DeclareMathOperator{\lcm}{lcm}
\DeclareMathOperator{\glcd}{glcd}

\DeclareMathOperator{\den}{den}
\DeclareMathOperator{\nr}{nr}
\DeclareMathOperator{\N}{N}
\DeclareMathOperator{\OC}{\mathrm{OC}}

\DeclareMathOperator{\SOC}{\mathrm{SOC}}
\DeclareMathOperator{\OS}{\mathrm{OS}}
\DeclareMathOperator{\SOS}{\mathrm{SOS}}

\begin{document}
\title{CSLs of the root lattice $\mathbf{A_4}$}

\author{M Heuer$^{1}$ and P Zeiner$^{2}$}

\address{$^{1}$ Department of Mathematics and Statistics, The Open University, Milton Keynes, UK}
\address{$^{2}$ Fakult{\"a}t f{\"u}r Mathematik, Universit{\"a}t Bielefeld, Germany}

\ead{m.heuer@open.ac.uk, pzeiner@math.uni-bielefeld.de}

\begin{abstract}
Recently, the group of coincidence isometries of the root lattice $A_4$ has been determined providing a classification of these isometries with respect to their coincidence indices. A more difficult task is the classification of all CSLs, since different coincidence isometries may generate the same CSL. In contrast to the typical examples in dimensions $d \leq 3$, where coincidence isometries generating the same CSL can only differ by a symmetry operation, the situation is more involved in $4$ dimensions. Here, we find coincidence isometries that are not related by a symmetry operation but nevertheless give rise to the same CSL.
We indicate how the classification of CSLs can be obtained by making use of the icosian ring and provide explicit criteria for two isometries to generate the same CSL. Moreover, we determine the number of CSLs of a given index and encapsulate the result in a Dirichlet series generating function.
\end{abstract}

\section{Introduction}
In crystallography, CSLs are a well established tool for the 
description and classification of grain boundaries, see \cite{baake} and references given there.
CSLs have been studied for several decades  in dimensions $d\leq 3$. In higher dimensions they became interesting after the discovery of 
quasicrystals. Here, we discuss a 4-dimensional lattice, namely the root lattice $A_4$, which is intimately 
related to quasicrystals with fivefold symmetry and the Penrose tiling in particular. 

Let us recall some basic notions first.
A \emph{lattice} $\G$ in the Euclidean space  $\RR^d$, can be described as the set of all integer linear combinations of a basis $\{b_1, \ldots, b_d\}$ of $\RR^d$, which is denoted by $\G=\langle b_1, \ldots, b_d \rangle_{\ZZ}$. A \emph{coincidence site lattice} (CSL) of a lattice $\G$ is defined as 
$\G  \cap \ts R\G$, where $R$ is a linear isometry and the corresponding \emph{coincidence index} 
$\varSigma (R) \, = \, [ \G : (\G\cap R\G)]\,$ is finite, i.e.\! $R$ and $R\G$ share a common sublattice. 
The orthogonal group, i.e. the group of all linear isometries of $\RR^d$, is denoted by  $\mathrm{O}(d,\RR)$. We define the set of all \emph{coincidence isometries} as
\[
    \OC (\G) \, := \,
    \{ R\in\mathrm{O}(d,\RR)\mid \varSigma (R) < \infty \}
\] and its restriction to rotations as $\SOC(\G)$.
According to \cite{baake, csl, GB}, $\OC(\G)$ and $\SOC(\G)$ are not only  subgroups of $\mathrm{O}(d,\RR)$, but   subgroups of the group of all similarity isometries which is defined as
$
    \OS (\G) \, := \, 
    \{ R\in\mathrm{O}(d,\RR)\mid \alpha R\G\subset\G \text{ for some }
    \alpha \in\RR^{+}\}\,. 
$  A sublattice of $\G$ of the form $\alpha R \G$ %, where $\alpha \in \RR^{+}$ and $R \in \OS(\G)$,
 is called a \emph{similar sublattice} (SSL) of $\G$.  
If we denote the
\emph{denominator} of a matrix $R\in\OS (\G)$ relative to the lattice
$\G$ as
$
   \den^{}_{\G} (R) \, := \,
   \min\ts \{\alpha \in  \RR^{+} \mid \alpha R\G\subset\G\}\, ,
$
we have the following characterisation, see \cite{csl} for details:

\begin{equation} \label{eq-OC in OS}
\OC (\G) = \{ R\in \OS (\G)\mid
  \den^{}_{\G} (R) \in \NN \}
\end{equation}

%When a lattice $\G \subset \RR^d$ is given, one is interested in the number of distinct SSLs of $\G$ of index $n$ as well as in the number of distinct CSLs of $\G$ of index $n$. If these arithmetic functions are multiplicative, they are often encapsulated into  Dirichlet series generating functions, because of their Euler product decomposition. A detailed introduction to arithmetic functions, the corresponding Dirichlet series and Euler products can be found in \cite{Apostol}. For many lattices in $d \leq 4$ the arithmetic functions which count the number of distinct SSLs and CSLs of each index have already been derived; see  for instance \cite{baake, BM98, BM99}. 
%One lattice for which this problem has not been solved completely is the root lattice $A_4$. In this paper we finally present the Dirichlet series which counts the number of distinct CSLs or the root lattice $A_4$ and its dual lattice $A^{*}_4$; details can be found in \cite{ssl, csl, BZ, Hphd}.

In the following sections we first describe the root lattice $A_4$ in a suitable setting. Then, we briefly introduce the key objects required for the analysis of its CSLs, as well as their relations to each other. We continue by reviewing some parts of the classification of the coincidence rotations of the root lattice $A_4$, see \cite{csl, H08} for details, and extend this analysis by providing explicit criteria for two coincidence isometries to  generate the same CSL. Finally, we classify the CSLs of the root lattice $A_4$ with respect to their coincidence indices and encapsulate the result in a Dirichlet series generating function which gives the number of CSLs for every index.

\section{The root lattice $A_4$}
The root lattice $A_4$ is usually defined in a $4$-dimensional hyperplane of $\RR^5$ as follows:
$$A_4:=\{(x_1, \ldots, x_5) \in \ZZ^5 \mid x_1+\ldots+x_5=0\}$$ 
However, this realisation of $A_4$ in $\RR^5$ is not suitable for our purposes. We prefer the following realisation in $\RR^4$
\begin{align*}\label{def-alat}
L:=\bigl\langle (1,0,0,0), \tfrac{1}{2}&(-1,1,1,1), (0,-1,0,0), \tfrac{1}{2}(0,1,\tau\!-\!1,-\tau) \bigr\rangle_{\ZZ}\,
\end{align*}
where $\tau=(1+\sqrt{5})/2$ is the golden ratio; see \cite{CMP,ssl, csl} for
details. Note that $L$ is just rescaled by a factor $\frac{1}{\sqrt{2}}$ in comparison to $A_4$. This particular description $L$ is very convenient for our problem, as it enables us to use the arithmetic of the quaternion algebra $\HH (\QQ(\sqrt{5}\,))$; see \cite{KR} for a detailed introduction to Hamilton's quaternions. To simplify the notation  we define 
\[
K:=\QQ(\sqrt{5}\,)=\{r+s \sqrt{5} \mid  r,s \in \QQ \},
\] which is a quadratic number field.
The algebra $\HH (K)$ is explicitly given as
\[
\HH (K) = K \oplus \ii  K \oplus \jj  K \oplus \kk  K,
\]
 where the generating elements satisfy Hamilton's relations 
$\ii^2 = \jj^2 =
\kk^2 = \ii\jj\kk = -1.$ 
It is equipped
with a \emph{conjugation} $\,\bar{.}\,$ which is the unique mapping that
fixes the elements of the centre of the algebra $K$ and reverses the
sign on its complement. If we write 
\[
q=(a,b,c,d)=a+\ii b+\jj c+\kk d, \quad \text{this means}\quad \bar{q}=(a,-b,-c,-d).
\]
The reduced norm and trace in $\HH (K)$ are defined by
\[   \nr (q) \, := \, q\bar{q} \, = \, \lvert q\rvert^2 
   \quad \text{and} \quad
   \tr (q) \, := \, q + \bar{q}, \, \]
where we canonically identify an element $\alpha\in K$ with the
quaternion $(\alpha,0,0,0)$.  For any $q\in\HH (K)$, $\lvert q \rvert$ is
its Euclidean length, which need not be an element of $K$.
Nevertheless, one has $\lvert rs\rvert = \lvert r\rvert \lvert
s\rvert$ for arbitrary $r,s\in\HH (K)$.  Due to the geometric meaning,
we use the notations $\lvert q\rvert^2$ and $\nr (q)$ in parallel. 
An element $q\in\HH (K)$ is called \emph{integral} when both $\nr (q)$
and $\tr (q)$ are elements of \[\oo:=\ZZ[\tau]:=\{m+n\tau \mid m, n \in \ZZ \},\] which is the ring of
integers of the quadratic number field $K$.

The \emph{icosian ring $\II$} consists of all linear combinations of the vectors 
\begin{align*} 
   (1,0,0,0),(0,1,0,0),
   \tfrac{1}{2}(1,1,1,1),\tfrac{1}{2}(1 \! - \!\tau,\tau,0,1)
   \end{align*} 
with coefficients in $\oo$. The elements of $\II$ are called \emph{icosians}. For more information on this remarkable object, see for example \cite{Reiner,MP,MW, ssl} and
references given there. 

An element $p\in\II$ is called \emph{$\II$-primitive} when
$\alpha p\in\II$, with $\alpha\in K$, is only possible with
$\alpha\in\oo$. Similarly, a sublattice $\L$ of $L$ is called \emph{$L$-primitive} when
$\alpha\L\subset L$, with $\alpha\in\QQ$, implies $\alpha\in\ZZ$. Whenever the context is clear, we simply use the term ``primitive'' in both cases.

The detailed arithmetic structure of $\II$ plays a key role in the characterisation of the coincidence rotations   for $L$; see \cite{csl}.
Another important object for this characterisation is the following map, called \emph{twist map}. If $q=(a,b,c,d)$, it is defined by the
mapping 
\[
q\mapsto \widetilde{q}:=\, (a{\ts}',b{\ts}',d{\ts}',c{\ts}')\, ,
\]
where ${}'$ denotes the algebraic conjugation in $K$, as defined by the
mapping $\sqrt{5}\mapsto -\sqrt{5}$. 
The algebraic conjugation in $K$ is
also needed to define the absolute norm on $K$, via 
$
     \N (\alpha) \, := \, \lvert \alpha \alpha{\ts}' \rvert\,  .
$
For the various properties of the twist map, e.g. $\ti{\ti{q}}=q$ and $\ti{pq}=\ti{q}\ti{p}$ where $p,q \in \II$, we refer the reader to
\cite{ssl} and references therein.  The most important properties in our
present context are summarised in the following Lemma, see \cite[Proposition 1]{ssl} and \cite[Lemma 4]{csl}, which describes the relations between $L$ and $\II$.
\begin{lemma} \label{fundamental}
Within $\HH (K)$, one has\/ $\widetilde{\II} = \II$ and
\[
L = \{ x\in\II\mid\widetilde{x}=x\}=\{ x + \widetilde{x}\mid x\in\II\} = \phi^{}_{+} (\II),
\] 
where the $\QQ$-linear
mapping $\phi^{}_{+} \! : \, \HH (K) \longrightarrow
\HH (K)$, is defined by $\phi^{}_{+} (x) = x+\widetilde{x}$.
\end{lemma}

The dual $A_4^{*}$ of the  the root lattice $A_4$, here in form of the dual of the lattice $L$, is given by
\[
L^{*}:=\{x \in \RR^{4} \mid \langle x|y \rangle \in \ZZ \text{ for all } y \in L\}.
\]

\section{CSLs}
The investigation of CSLs of $L$  can be restricted to rotations only, since $\overline{L}=L$, which means that any orientation reversing operation can be obtained from an orientation preserving one after applying conjugation first. Let $\SOS(L)$ be the group of all similarity rotations of the lattice $L.$
Considering Eq. \eqref{eq-OC in OS}, it is obvious how 
$\SOC (L)$ and $\SOS (L)$ are related in general. 
By \cite[Corollary 1]{ssl} we know that any
similarity rotation of $L$ is of the form $x\mapsto
\frac{1}{\lvert q\ti{q}\ts\rvert}\, q x \ti{q}$ with
$q\in\II$, which we refer to as $R(q)$, or
more precisely we define $R(q)x:= \frac{1}{\lvert q \ti{q}\ts\rvert}\,q x \ti{q}$. We denote the denominator of $R(q)$ by $\den(q)$. Among the elements of $\SOS(L)$, we can identify the elements of $\SOC (L)$ as follows; see \cite{csl} for the proof.

\begin{proposition} \label{proposition-commensurate iff den(q) integer}
  Let $0\neq q\in\II$ be an arbitrary icosian.  Then, the lattice
  $\frac{1}{\lvert q\ti{q}\ts\rvert}\, q L \ti{q} \cap L$
  is a CSL of $L$ if and only if\/ $\lvert q
  \ti{q}\ts\rvert\in\NN$. If $q$ is primitive, then $\den(q)=|q \ti{q}|$.
\end{proposition}
We call an icosian $q\in\II$ \emph{admissible} when
$\lvert q \ti{q}\ts\rvert\in\NN$. As $\nr (\ti{q}\ts) =
\nr (q){\ts}'$, the admissibility of $q$ implies that
$\N\bigl(\nr (q)\bigr)=\lvert q \ti{q}\ts\rvert^2$ is a square in $\NN$. 
An immediate consequence of Proposition \ref{proposition-commensurate iff den(q) integer} is that
\begin{equation} \label{eq-parametrisation of SOC(L)}
 \SOC(L)= \set{R(p) \mid p \in \II \text{ is primitive and admissible }},
\end{equation}
which delivers the following classification.

\begin{theorem} \label{theorem-Form of CSL of L}
  The CSLs of $L$ are precisely the lattices of the form $L
    \cap \frac{1}{\lvert q\ti{q}\ts\rvert}\, qL\ti{q}$
    with $q\in\II$ primitive and admissible.
\end{theorem}

This is the first step to connect certain primitive right
ideals $p\II$ of the icosian ring with the CSLs of $L$.
Before we continue in this direction, let us consider the
relation with the coincidence rotations, see \cite[Lemma 5]{csl}.

\begin{lemma} \label{symmetries} 
  Let $r,s\in\II$ be primitive and admissible quaternions, with
  $r\II=s\II$.  Then, one has $L\cap \frac{rL\ts
    \tilde{r}}{\lvert r \tilde{r}\rvert} = L\cap \frac{sL
    \tilde{s}}{\lvert s \tilde{s}\rvert}$.
\end{lemma}

The converse statement to
  Lemma~\ref{symmetries} is not true, as the equality of
  two CSLs does \emph{not} imply the corresponding ideals to be equal. An example is provided by $r=(\tau,2\tau,0,0)$ and
  $s=(\tau^2,\tau,\tau,1)$, which define the same CSL, though
  $s^{-1}r$ is not a unit in $\II$. The CSL is spanned by the
  basis $\{(1,2,0,0),(2,-1,0,0),(\frac{3}{2},\frac{1}{2},\frac{1}{2},
  \frac{1}{2}),(-1,\frac{1}{2},\frac{\tau-1}{2},-\frac{\tau}{2})\}$.
  Note that, as a direct
  consequence of  \cite[Lemma 5]{ssl}, two primitive 
  quaternions $r,s\in\II$ are related by a
  rotation symmetry of $L$ if and only if $r\II=s\II$.

For our further discussion we need to replace the primitive and admissible icosian $p$, and with it $\tilde{p}$, by certain $\oo$-multiples, such that their norms have the
same prime divisors in $\oo$. In view of the form of the rotation $x\mapsto
\frac{1}{\lvert p\widetilde{p}\ts\rvert}\, p x \widetilde{p}$, this is actually rather natural because it
restores some kind of symmetry of the expressions in relation to
the two quaternions involved. For a primitive and admissible icosian $p\in\II$  we
choose explicitly 
\begin{align} \label{def-alpha}
    \alpha^{}_{q} \, = \, 
    \sqrt{\frac{\lcm(\nr (q),\nr (\widetilde{q}\,))}{\nr (q)}}
    \, \in \, \oo \, ,
\end{align}
where we assume a suitable standardisation for the $\lcm$ (least common multiple) of
two elements of $\oo$; see \cite{csl} for details. Moreover, we have the relation $\alpha^{}_{\tilde{q}}
= \widetilde{\alpha^{}_{q}} = \alpha^{\,\prime}_{q}$.
The icosian $\alpha_{q} q$ is called the \emph{extension} of the
primitive admissible element $q\in\II$, and $(\alpha^{}_q q,
\alpha^{\,\prime}_{q} \widetilde{q}\,)$ the corresponding
\emph{extension pair}.  Since $\alpha_q$ and $\alpha_q'$ are central, the extension does not change the rotation, i.e.
\begin{align*} 
     \frac{q x \widetilde{q}}{\lvert q\widetilde{q}\ts\rvert} 
     \, = \, \frac{q_{\alpha} x \widetilde{q_{\alpha}}}
     {\lvert q_{\alpha}\widetilde{q_{\alpha}}\ts\rvert}
\end{align*}
holds for all quaternions $x$.  Note that the definition of the
extension pair is unique up to units of $\oo$, and that one has the
relation
\begin{align*} 
     \nr (q_{\alpha}) \, = \,
     \lcm \bigl(\nr (q), \nr(\widetilde{q}\,)\bigr)
     \, = \, \nr (\widetilde{q_{\alpha}})
     \, = \, \lvert q_{\alpha}\ts
             \widetilde{q_{\alpha}}\rvert
     \, \in \, \NN \, .
\end{align*}
Following \cite{csl} we define the set
\[
     L (q) \, = \, \set{ qx + \ti{x}\ti{q}
     \mid x\in\II } \, = \, \phi^{}_{+} (q\ts\II)\, .
\]
Note
that, due to $\ti{\II}=\II$, one has $L (q) =
\ti{L (q)}$.
This leads to the following parametrisation of the CSLs of $L$, compare \cite[Theorem 2]{csl}. 
\begin{theorem} \label{theorem-CSL of L equal L(q)}
   Let $q\in\II$ be admissible and primitive,
   and let $q_{\alpha} = \alpha_{q}\ts q$ be its extension.
   Then, the CSL defined by $\ts q$ is given by
\[
   L \cap \frac{1}{\lvert q\ti{q}\ts\rvert}\ts
    q L\ti{q} \, = \, L (q_{\alpha})\, = \,(q_{\alpha}\II + \II \ti{q}_{\alpha}) \cap L.
\]
and its index is 
$
  \Sigma_{L}(q) \, = \, \nr(q^{}_{\alpha}) \, = \,
   \lcm (\nr(q),\nr(q)') \, .
$
\end{theorem}

\section{The Number of CSLs of each Index}
Until now we have followed the argumentation in \cite{csl}, which lead to the classification of all coincidence rotations of the lattice $L$. In order to obtain a classification of the actual CSLs, we make use of the corresponding classification for the icosian ring $\II$ and its submodul $L[\tau]$, which can  either be seen as  $\oo$-modules of rank $4$ or as lattices in $\RR^8$, so all our definitions of the introductory section apply. Since any coincidence rotation of $L$ is a coincidence rotation of $L[\tau]$ and $\II$, we can  start our analysis from the classification of the CSLs of $\II$, which directly gives a first sufficient condition when two CSLs of $L$ are equal, see \cite{Hphd} for details.

\begin{lemma} \label{lemma-If nr and glcd are equal CSL of L are equal}
 Let $p_1, p_2 \in \II$ be primitive and admissible, such that 
\[|p_1|^2=|p_2|^2  \and \glcd(p_1, |p_1\ti{p}_1|)=\glcd(p_2, |p_2\ti{p}_2|).
\] Then $L \cap \frac{1}{\rvert p_1 \ti{p_1}\lvert} p_1 L \ti{p_1}=L \cap \frac{1}{\rvert p_2 \ti{p_2}\lvert} p_2 L \ti{p_2}.$
\end{lemma}
In fact, the converse holds only if $|p_1|^2=|p_2|^2$ is not divisible by
the splitting prime $p=5$. For the general situation we have the following

\begin{theorem}
 Let $p_1, p_2 \in \II$ be primitive and admissible. 
Then \[ 
L \cap \tfrac{1}{\rvert p_1 \ti{p_1}\lvert} p_1 L \ti{p_1}=L \cap \tfrac{1}{\rvert p_2 \ti{p_2}\lvert} p_2 L \ti{p_2}
\]
if and only if 
\begin{equation}\label{eq-condition}
|p_1|^2=|p_2|^2  \and \glcd(p_1, \tfrac{|p_1\ti{p}_1|}{c})=\glcd(p_2, \tfrac{|p_2\ti{p}_2|}{c}),
\end{equation}
where $c=\sqrt{5}$, if $|p_1|^2=|p_2|^2$ is divisible by $5$ and $c=1$ otherwise.
\end{theorem}

Note that $|p_1|^2=|p_2|^2$ directly implies that the coincidence indices and denominators coincide, which is clearly  necessary for the CSLs to be equal. The second condition of \eqref{eq-condition} tells us how much $p_1$ and $p_2$ are allowed to differ in terms of their decomposition in irreducible elements of $\II$. If
$|p_1|^2=|p_2|^2$ is not divisible by $5$ and $\nr(p_1)=\nr(p_1)'$, then
(and only then) the second condition is equivalent to $p_1=p_2\epsilon$ where $\epsilon$ is a unit of $\II$, which means that $R_1$ and $R_2$ are symmetry related. Hence, in this special case, and only then, we can deduce from two equal CSLs that they are generated by coincidence rotations which are symmetry related, whereas in all other cases there are at least two rotations that
are not symmetry related but give the same CSL.

Now, we can conclude with the derivation of the Dirichlet series for the number of CSLs of $L$ of a given index, see \cite{apostol} for a introduction to Dirichlet series. Note, that by \cite[Lemma 2.5, Theorem 2.2]{baake} the number of CSLs of a given index $n$ is the same for the lattices $A_4$ and $L$ as well as their dual lattice $L^*$. Let $f(n)$ be the number of CSLs of
  index $n$. Then, $f(n)$ is a multiplicative
  arithmetic function, given at prime powers $p^r$ with $r \geq 1$ by

\[
f(p^r)=
\begin{cases}
6 \cdot 5^{2r-2} & \mbox{ for } p=5 \\
(p^{2}+1)p^{2r-2} & \mbox{ for } p\equiv \pm 2 \pmod 5 \\
\frac{(p+1)^2}{p^3-1}(p^{2r+1}+p^{2r-2}-2p^{\frac{(r-1)}{2}})
& \mbox{ for } p\equiv \pm 1 \pmod 5, r \mbox{ odd}, \\
\frac{(p+1)^2}{p^3-1}
\left(p^{2r+1}+p^{2r-2}-2\frac{(p^2+1)}{p+1}p^{\frac{(r-2)}{2}}\right)
& \mbox{ for } p\equiv \pm 1 \pmod 5, r \mbox{ even}.
\end{cases}
\]
Finally elementary calculations involving geometric series yield its Dirichlet series generating function compare 
\cite{BZ, Hphd}:
\[
\begin{aligned} 
D (s)  & = \sum_{n=1}^{\infty}\frac{f(n)}{n^s} 
 = {\textstyle \,\left(1+ 6 \frac{ 5^{-s}}{1-5^{2-s}}\, \right)} \prod_{p\equiv\pm 1 \,(5)}
    {\textstyle \frac{1+p^{-s}+2p^{1-s}+2p^{-2s}+p^{1-2s}+p^{1-3s}}{(1-p^{2-s})\,(1-p^{1-2s})}}\,
    \prod_{p\equiv\pm 2 \,(5)} {\textstyle \frac{1+p^{-s}}{1-p^{2-s}}}  \\
  &= \,{\textstyle 1 + \frac{5}{2^s} + \frac{10}{3^s} + \frac{20}{4^s} + 
      \frac{6}{5^s} + \frac{50}{6^s} + \frac{50}{7^s} + \frac{80}{8^s}+ \frac{90}{9^s} +  \frac{30}{10^s} + \frac{144}{11^s} +  \ldots}
\end{aligned}
\]

\ack{Acknowledgements}
It is a pleasure to thank M. Baake and U. Grimm for their suggestions and helpful discussions. This work was supported by EPSRC grant EP/D058465/1 and the German Research Council (DFG), within the CRC 701.
\bibliography{references}

\end{document}